\documentclass[12pt]{amsart}
\usepackage[margin=1in,letterpaper]{geometry}
\usepackage[T1]{fontenc}
\usepackage{lmodern,amsmath,amssymb,amsfonts,amsthm,mathrsfs}
\usepackage{mathtools,booktabs,microtype,needspace, color}
\usepackage[shortlabels]{enumitem}
\setlist[enumerate]{itemsep=1em,parsep=0pt,topsep=0.5em}
\newlength{\savedenumitemsep}
\NewDocumentCommand{\itemafterdisplay}{o}{%
  \setlength{\savedenumitemsep}{\itemsep}%
  \setlength{\itemsep}{0pt}%
  \IfNoValueTF{#1}{\item}{\item[#1]}%
  \setlength{\itemsep}{\savedenumitemsep}%
}
\usepackage[hidelinks]{hyperref}
\hypersetup{pdftitle={Nonnegative weak tilings and the Coven--Meyerowitz condition (T2)},pdfauthor={Gergely Kiss, Caleb Marshall, M\'at\'e Matolcsi, G\'abor Somlai}}
\newtheorem{theorem}{Theorem}[section]
\newtheorem{proposition}[theorem]{Proposition}
\newtheorem{lemma}[theorem]{Lemma}
\newtheorem{corollary}[theorem]{Corollary}
\theoremstyle{definition}
\newtheorem{definition}[theorem]{Definition}
\newtheorem{question}[theorem]{Question}

\newtheorem{remark}[theorem]{Remark}
\numberwithin{equation}{section}
\newcommand{\zz}{\mathbb Z}
\newcommand{\qq}{\mathbb Q}
\newcommand{\rr}{\mathbb R}
\newcommand{\cc}{\mathbb C}

\newcommand{\one}{\mathbf 1}
\newcommand{\calm}{\mathcal M}
\newcommand{\supp}{\operatorname{supp}}

\newcommand{\mass}{\operatorname{mass}}
\newcommand{\CRT}{\operatorname{CRT}}
\newcommand{\Ph}{\Phi}

\title[Nonnegative weak tilings and (T2)]{Nonnegative weak tilings and the Coven--Meyerowitz condition (T2)}
\author{Gergely Kiss}
\address{Gergely Kiss \newline
HUN-REN Alfr\'ed R\'enyi Institute of Mathematics, Analysis Department, Reáltanoda u. 13-15, 1053, Budapest, Hungary, \newline
and \newline
Corvinus University of Budapest, F\H{o}v\'am t\'er 13-15, 1093, Budapest, Hungary}
\email{kigergo57@gmail.com}
\author{Caleb Marshall}
\address{Caleb Marshall, \newline University of Toronto, Department of Mathematics, St. George Campus, Toronto, ON, Canada}
\email{caleb.marshall@utoronto.ca}
\author{M\'at\'e Matolcsi}
\address{M\'at\'e Matolcsi, \newline HUN-REN Alfréd Rényi Institute of Mathematics, Reáltanoda u. 13-15, 1053, Budapest, Hungary,
\newline
and \newline Budapest University of Technology and Economics, Department of Analysis and Operations Research,
Institute of Mathematics,
M\H{u}egyetem rkp. 3., H-1111 Budapest, Hungary.}
\email{matolcsi.mate@renyi.hu}
\author{G\'abor Somlai}
\address{G\'abor Somlai,
The University of Melbourne, School of Mathematics and Statistics, Parkville VIC 3010, Australia,
\newline
and \newline
on unpaid leave from E\"otv\"os Lor\'and University}
\email{gabor.somlai@unimelb.edu.au}
\date{October 20, 2026}

\begin{document}
\begin{abstract}
We develop a general method, in finite cyclic groups, for constructing nonnegative rational functions $f$, normalized by $f(0)=1$, that satisfy
the Coven--Meyerowitz condition $(T1)$ but fail $(T2)$, such that $f$ is a weak
complement for a set tile $A$. Motivated by the unsupported-divisor questions
in~\cite{KLMS2026a}, we show that if a set tile $A$ admits an unsupported
divisor, then a suitable dilation of $A$ admits such a rational weak
complement $f$. Using this method, we construct an explicit example in the group $\zz_{180}$. The weak complements in our
explicit perturbative constructions are not positive definite because of
restrictions on the Fourier support of the perturbations. Thus these
constructions do not resolve the proposed extension of the
Coven--Meyerowitz conditions to positive-definite weak
tilings~\cite[Problem~1.2]{KLMS2025}.
\end{abstract}
\maketitle

\section{Introduction}\label{sec:intro}

By D.~Newman's periodicity theorem~\cite{N1977}, tilings of the integers
by a finite set can be studied through factorizations
$A\oplus B=\zz_M$, in which each residue has a unique representation
as $a+b$. For $A\subset\zz_M$, choose representatives in
$\{0,\ldots,M-1\}$ and write $A(X)=\sum_{a\in A}X^a$ for its
\emph{mask polynomial}. The tiling identity becomes
\begin{equation}\label{eq:introtiling}
 A(X)B(X)\equiv U_M(X)\mod (X^M-1),\qquad
 U_M(X)=1+X+\cdots+X^{M-1}.
\end{equation}

E.~Coven and A.~Meyerowitz~\cite{CM1999} introduced two conditions on the
cyclotomic divisors of mask polynomials. Condition $(T1)$ determines
the value at $X=1$ from the prime-power divisors; $(T2)$ requires
products of prime-power orders associated with distinct primes to
occur as well. Together they imply tiling. Every integer tile
satisfies $(T1)$, and $(T2)$ is also necessary when its cardinality
has at most two distinct prime divisors (see \cite[Theorems A, B1, and B2]{CM1999}). Precise definitions of these properties are
given in Section~\ref{sec:background} of this work.

\begin{remark}
Recent preprints by J.~Prajapati~\cite{P2026} and
C.~Li et al.~\cite{LWHWWZ2026} claim the unrestricted necessity of $(T2)$
for integer tiles, which, if correct, would resolve the Coven--Meyerowitz
conjecture for integer tilings. Our work does not use these proposed
solutions or their methods.
\end{remark}

Our starting point is the study of unsupported cyclotomic divisors in
G.~Kiss, I.~{\L}aba, C.~Marshall, and G.~Somlai~\cite[Section~1.3]{KLMS2026a}. They
construct tiles satisfying $(T1)$ and $(T2)$ with unsupported
divisors~\cite[Theorem~1.11]{KLMS2026a}. We use such mixed zeros to
construct rational weak complements failing $(T2)$.

Following~\cite{KLM2023,LM2022}, a \emph{normalized nonnegative weak complement}
of $A\subset\zz_M$, with $0\in A$, is a function satisfying
\begin{equation}\label{eq:weakintro}
 f:\zz_M\longrightarrow\rr_{\ge0},\qquad
 f(0)=1,\qquad \one_A*f=\one_{\zz_M}.
\end{equation}
Equivalently, $\one_A*(f-\delta_0)=\one_{A^c}$, so $A$ weakly tiles
its complement. We keep the set factor equal to a tile and allow
rational weights in its complement, leading to the following question.

\begin{question}\label{q:weak}
Suppose that $A\subset\zz_M$ is a tile satisfying $(T1)$ and $(T2)$.
Must every normalized nonnegative rational weak complement of $A$ that
satisfies $(T1)$ also satisfy $(T2)$?
\end{question}

The answer is negative, as we shall see in the sequel.
However, the explicit
weak complements in our construction are not positive definite; the distinction from positive-definite (PD) weak
tilings is discussed in Section~\ref{sec:pdlimitations}.

\subsection{An explicit counterexample}
We begin with an example of period $180$. The polynomial appearing
below is given explicitly in Section~\ref{sec:example180}.

\begin{samepage}
\begin{theorem}\label{thm:introexample}
The set
\[
 A_*=\{0,12,18,24,36\}\subset\zz_{180}
\]
tiles $\zz_{180}$ and satisfies ordinary $(T1)$ and $(T2)$. There is
an ordinary polynomial $F\in\qq[X]$ with nonnegative coefficients and
$F(0)=1$ such that
\[
 A_*(X)F(X)\equiv U_{180}(X)\mod (X^{180}-1),
\]
while $F$ satisfies ordinary $(T1)$ but fails $(T2)$.
Its cyclic coefficient function
$w_F^{180}:\zz_{180}\longrightarrow\qq_{\ge0}$, which satisfies
\[
 F(X)\equiv\sum_{y=0}^{179}w_F^{180}(y)X^y\mod (X^{180}-1),
\]
is not positive definite.
\end{theorem}
\end{samepage}

The set $A_*$ represents each residue class modulo five exactly once,
so that the subgroup $B_*=5\zz_{180}$ is a tiling complement of
cardinality $36$. To obtain $F$, we use a different complement, formed
from translated layers, and verify the required cyclotomic divisibility
properties directly in Section~\ref{sec:example180}. We do not claim
that $180$ is the least possible period.

The example comes from the five-point set
$E=\{0,2,3,4,6\}$, since $A_*=6E$ and
$E(X)=\Ph_5(X)\Ph_6(X)$. This is an instance of the recombination
construction in~\cite[Example~8.1]{KLMS2026a}, with primes $2,3,5$.
H.~Tan and Y.~Zhang~\cite{TZ2026} also use this set $E$ in their
unsupported-divisor constructions.

\subsection{A construction without new primes}
A \emph{mixed order} is an integer with at least two distinct prime
divisors. If $N=\prod_i p_i^{\alpha_i}$, its
\emph{constituent prime-power orders} are the $p_i^{\alpha_i}$.
The following theorem gives a general construction from a mixed
cyclotomic divisor whose constituent prime-power divisors are absent.

\Needspace{12\baselineskip}
\begin{samepage}
\begin{theorem}\label{thm:introconstruction}
Suppose that $A\oplus B=\zz_M$ and that $N\mid M$ is a mixed order.
Assume that $\Ph_N\mid A(X)$, but that $\Ph_s\nmid A(X)$ for every
constituent prime-power order $s$ of $N$.

Then there exist an integer $q\ge1$ and a polynomial $P\in\qq[X]$ with
nonnegative coefficients and $P(0)=1$ such that
\[
 A(X^q)P(X)\equiv U_{qM}(X)\mod (X^{qM}-1),
\]
and $P$ satisfies ordinary $(T1)$ but fails $(T2)$.
The dilation $q$ can be chosen to be $1$ or a power of a prime dividing
$N$.
\end{theorem}
\end{samepage}

A more thorough description of the polynomial $P$ and the dilated set tile $qA \subset \mathbb{Z}_{qM}$ with mask polynomial $A(X^q)$ is given in Theorem~\ref{thm:primepreserving}.

The tiles constructed by G.~Kiss, I.~{\L}aba, C.~Marshall, and
G.~Somlai~\cite[Theorem~10.1]{KLMS2026a} satisfy the hypotheses of
Theorem~\ref{thm:introconstruction}, so their dilations admit the weak
complements constructed here. The set $qA\subset\zz_{qM}$ is again a set
tile. Furthermore, $qA$ satisfies ordinary $(T2)$ if and only if $A$
does, as shown in Theorem~\ref{thm:primepreserving}; see
also~\cite[Lemma~1.4]{CM1999}. Thus dilation preserves the cardinality,
tiling property, and Coven--Meyerowitz conditions of the set factor.

For the main step of the construction, we choose $q>1$ and form a set
complement $C$ of $qA$ by placing translated copies of $B$ in the
residue classes modulo $q$. Its \emph{normalized autocorrelation},
$|C|^{-1}\one_C*\one_{-C}$, is a weak complement of $qA$. We then
add a signed function whose Fourier support consists of a single
mixed order. This function, which we call a \emph{pure-order cuboid
filter}, is defined in Section~\ref{sec:filter}. The translations are
chosen so that the autocorrelation has enough positive weight to
absorb the negative coefficients of the filter, while the perturbation
does not change the value at the origin. The proof is given in
Section~\ref{sec:cycliclift}. The example in
Theorem~\ref{thm:introexample} uses the same idea with the composite
dilation $q=6$ and fewer layers than the general bound
would give.

Our use of a larger ambient group was inspired by the stacking
construction of R.~Greenfeld and T.~Tao~\cite[Theorem~3.1]{GT2024}.
As part of their disproof of the periodic tiling conjecture, they use
additional layers to encode a system of tiling equations as a single
equation in a larger group. Here the additional residue layers allow us
to place translated complements so that their autocorrelation absorbs
the signed perturbation while preserving normalization.

\subsection{The weak complement construction can fail without dilations}
One may ask whether the dilation in
Theorem~\ref{thm:introconstruction} can be avoided by allowing a larger
ambient period while leaving the set tile $A$ fixed. The following result
shows that this is not possible in general.

\begin{samepage}
\begin{theorem}\label{thm:introobstruction}
Let
\[
 T=\{0,2,3,4,6\}+30\{0,1,2,3,4,5\}\subset\zz_{180}.
\]
This thirty-element tile satisfies ordinary $(T1)$ and $(T2)$ and has
unsupported divisor $\Ph_6$.

For every $L\ge2$, suppose that $P\in\qq[X]$ has nonnegative
coefficients, satisfies $P(0)=1$ and ordinary $(T1)$, and obeys
\[
 T(X)P(X)\equiv U_L(X)\mod (X^L-1).
\]
Then $P$ satisfies ordinary $(T2)$.
\end{theorem}
\end{samepage}

Here \emph{unsupported divisor} includes the cardinality condition in
Definition~\ref{def:unsupported}; Theorem~\ref{thm:introconstruction} does not require that
additional condition. Theorem~\ref{thm:introobstruction} is proved as
Proposition~\ref{prop:fixedmask}.

We also consider what dilation can accomplish when the complement is
required to be a set. Proposition~\ref{prop:dilationequivalence}
shows that, for every tile $A\subset\zz_M$ and every $q\ge1$, the
set $qA\subset\zz_{qM}$ has a set complement failing $(T2)$ if and
only if $A$ has one. Thus dilation cannot produce a new failure at the
level of set tilings.

There is a further restriction on the ambient period. If $M$ has at
most two distinct prime divisors and $A\subset\zz_M$ is a tile, then
every rational polynomial weak complement of $A$ satisfying ordinary
$(T1)$ also satisfies ordinary $(T2)$
(Proposition~\ref{prop:twoprime}). In fact, this conclusion requires
neither nonnegativity nor normalization of the complement.

Finally, we explain why the perturbation method does not give
positive-definite examples. In the nontrivial constructions, the
perturbation vanishes at the origin and its Fourier transform is
supported on the Fourier zeros of the initial PD complement. Under these
assumptions, Proposition~\ref{prop:pdperturbation} shows that every
nonzero perturbation destroys positive definiteness. For the five-point
set underlying Theorem~\ref{thm:introexample}, we prove a stronger
statement: for every $q\ge1$, a normalized nonnegative rational PD
complement of $qE\subset\zz_{30q}$ satisfying ordinary $(T1)$ must
also satisfy ordinary $(T2)$ (Corollary~\ref{cor:pdset}).

\subsection*{Organization}
Section~\ref{sec:background} gives the Fourier and cyclotomic
conventions. Section~\ref{sec:filter} develops the cuboid filter, and
Sections~\ref{sec:cycliclift} and~\ref{sec:example180} prove Theorems~\ref{thm:introconstruction}
and~\ref{thm:introexample}. Section~\ref{sec:obstructions} proves Theorem~\ref{thm:introobstruction} and the
additional boundary results.
Section~\ref{sec:pdlimitations} explains the PD obstruction.
Appendix~\ref{sec:lift} presents an alternative construction using an
auxiliary prime coprime to the original period.

\section{Mask polynomials, Fourier transforms, and tiling conditions}\label{sec:background}

\subsection{Cyclic coefficient functions and convolution}
We will use both polynomials and their coefficient functions on finite
cyclic groups, extending the multiset notation
of~\cite[Section~2.1]{KLMS2026a} to rational weights. To fix the
correspondence, write $[X^n]Q$ for the
coefficient of $X^n$ in a polynomial $Q$. For $Q\in\qq[X]$ and $M\ge2$,
define its \emph{cyclic coefficient function} by
\[
 w_Q^M(x)=\sum_{\substack{n\ge0\\n\equiv x\bmod M}}[X^n]Q,
 \qquad x\in\zz_M.
\]
Thus coefficients whose exponents agree modulo the period are added
together, and the canonical representative of $Q$ modulo $X^M-1$ is
$\sum_{x=0}^{M-1}w_Q^M(x)X^x$. We call such a weight function a rational
multiset on $\zz_M$, and denote the collection by $\calm_{\qq}(\zz_M)$,
using $\calm_{\qq}^+(\zz_M)$ when all weights are nonnegative. Sets
correspond to weights in $\{0,1\}$.

Conversely, for a function $f:\zz_M\to\qq$, its \emph{canonical mask polynomial} $F \in \mathbb{Q}[X]$ is given by the formula
\[
F(X)=\sum_{x=0}^{M-1}f(x)X^x.
\]
Thus $F(0)=f(0)$ and
$F(1)=\sum_x f(x)$; polynomial evaluation at one is distinct from
the value $f(1)$ of the coefficient function.

For functions $f,g:\zz_M\to\cc$, write
\[
 (f*g)(x)=\sum_{y\in\zz_M}f(y)g(x-y),\qquad
 \mass(f)=\sum_x f(x),\qquad \widetilde f(x)=\overline{f(-x)}.
\]
Here the overline denotes complex conjugation. With these conventions,
multiplication of mask polynomials corresponds to convolution over the corresponding cyclic group, while reduction
modulo $X^M-1$ preserves mass: $Q(1)=\mass(w_Q^M)$.

\begin{samepage}
\begin{lemma}\label{lem:integer}
If $A\subset\zz_M$ is nonempty and $b:\zz_M\to\zz_{\ge0}$ satisfies
$\one_A*b=\one_{\zz_M}$, then $b$ is the indicator function of a set.
\end{lemma}
\end{samepage}
\begin{proof}
Fix $a_0\in A$. Then, for every $x$, one must have
$$
1 = \big(\mathbf{1}_A * b\big) (x+a_0) \coloneqq \sum\limits_{a \in A} b(x + a_0 - a) \geq b(x) \geq 0.
$$
Hence, $b(x) \in \{0,1\}$ for every $x \in \mathbb{Z}_M$, as claimed.
\end{proof}

The lemma explains why passing from sets to nonnegative integer weights
does not enlarge the class of complements at tiling level one. For a
rational weak complement, however, clearing a denominator produces a
\emph{multiple} tiling: if $qf$ is integer-valued, then
$\one_A*(qf)=q\one_{\zz_M}$.

\subsection{Fourier and cyclotomic criteria}
The tiling congruence can be tested at roots of unity, or equivalently
on the Fourier side. We use the unnormalized Fourier transform:
\[
 \widehat f(k)=\sum_{x=0}^{M-1}f(x)e^{-2\pi ikx/M},\qquad
 f(x)=\frac1M\sum_{k=0}^{M-1}\widehat f(k)e^{2\pi ikx/M}.
\]
The character indexed by $k$ has order $M/\gcd(k,M)$, with the trivial
character having order one. Character orthogonality gives the inversion
formula and the identity $\widehat{f*g}=\widehat f\widehat g$. To express
Fourier vanishing in terms of polynomial divisibility, recall that for
$d\ge1$ the cyclotomic polynomial is
\[
 \Ph_d(X)=\prod_{\substack{1\le u\le d\\(u,d)=1}}
       (X-e^{2\pi iu/d}).
\]
It is monic and irreducible over $\qq$, and
$X^M-1=\prod_{d\mid M}\Ph_d(X)$.

\begin{samepage}
\begin{proposition}[cf.~E.~Coven and A.~Meyerowitz, {\cite[Lemma~1.3]{CM1999}}]\label{prop:criterion}
For $A\subset\zz_M$ and a rational-valued function $f$ with canonical mask polynomial $F$,
the following are equivalent:
\begin{enumerate}[(i)]
\item $\one_A*f=\one_{\zz_M}$;
\item $A(X)F(X)\equiv U_M(X)\mod (X^M-1)$;
\item $|A|F(1)=M$ and, for every $d\mid M$ with $d>1$, at least one of
$A(X)$ and $F(X)$ is divisible by $\Ph_d(X)$.
\end{enumerate}
\end{proposition}
\end{samepage}
\begin{proof}
Comparing cyclic coefficients gives the equivalence of (i) and (ii).
On the Fourier side, the constant function has value $M$ at the trivial
character and zero elsewhere. The convolution identity is therefore
equivalent to the mass identity together with
$\widehat{\one_A}(k)\widehat f(k)=0$ at every nontrivial character.
Since the mask polynomials have rational coefficients, vanishing at one primitive
$d$th root is equivalent to divisibility by $\Ph_d$. This proves (iii).
\end{proof}

\subsection{The two conventions for \texorpdfstring{$(T1)$ and $(T2)$}{(T1) and (T2)}}
For rational mask polynomials it is useful to distinguish the conditions on an
ordinary polynomial from the conditions seen at a fixed cyclic period.
The latter involve only orders dividing that period.
\begin{definition}[E. Coven and A. Meyerowitz {\cite{CM1999}}]\label{def:CM}
For a nonzero polynomial $Q\in\qq[X]$, let
\[
 S_Q^*=\{p^a:\ p\text{ prime},\ a\ge1,\ \Ph_{p^a}\mid Q\}.
\]
The ordinary Coven--Meyerowitz conditions are:
\begin{enumerate}[(T1)]
\item $Q(1)=\prod_{s\in S_Q^*}\Ph_s(1)$;
\item whenever $k\ge1$ and $s_1,\ldots,s_k\in S_Q^*$ are powers of distinct primes,
$\Ph_{s_1\cdots s_k}\mid Q$.
\end{enumerate}
For a function on $\zz_M$, we call the versions restricted to the ambient
period the \emph{period-relative} conditions. They use
$S_{Q,M}=\{s\in S_Q^*:s\mid M\}$ in the same formulas instead of $S_Q^*$. These depend only on the cyclic
coefficient function, not on its ordinary polynomial representative.
\end{definition}

Recall that $\Ph_{p^a}(1)=p$.
We write $p^n\Vert M$ when $p^n\mid M$ and $p^{n+1}\nmid M$.
For a prime $p$ and a positive integer $m$, the \emph{$p$-adic valuation}
$v_p(m)$ is the unique integer $r\ge0$ such that $p^r\Vert m$. The following explains how prime-power divisors $p^{\alpha} \mid M$ of the tiling period $M$ are distributed among the mask polynomials of the sets.

\begin{samepage}
\begin{lemma}\label{lem:allocation}
If $A\oplus B=\zz_M$, then every prime-power cyclotomic divisor of either
mask polynomial has order dividing $M$. For each prime power $p^a\mid M$,
$a\ge1$, the factor $\Ph_{p^a}$ divides exactly one mask polynomial.
Both mask polynomials satisfy ordinary $(T1)$.
\end{lemma}
\end{samepage}
\begin{proof}
The counting argument rests on integrality. For a nonzero $Q\in\zz[X]$
with $Q(1)>0$, monicity and Gauss's lemma imply that the product of any
distinct cyclotomic divisors divides $Q$ in $\zz[X]$. Evaluating at one
shows that the number of prime-power factors based at $p$ is at most
$v_p(Q(1))$.

Now let $p^n\Vert M$. Proposition~\ref{prop:criterion} requires each of
the $n$ factors $\Ph_p,\ldots,\Ph_{p^n}$ to divide at least one of the
two masks. On the other hand, their cardinalities allow only
$v_p(|A||B|)=n$ such factors in total. Thus every required factor
occurs exactly once, and there is no room for any further $p$-power
factor. A prime not dividing $M$ cannot occur in either mask polynomial for the
same reason. Evaluating the assigned factors at one gives $(T1)$.
\end{proof}

\begin{remark}\label{rmk:whyT1}
The restriction to integer mask polynomials is essential to this counting
argument, as a small example shows. The tile
$\{0,3\}\subset\zz_6$ has the normalized weak complement
$f=(1,3/4,1/4,0,1/4,3/4)$. Its mass is $3$, but it vanishes at orders
$2$ and $3$, whose cyclotomic values at one multiply to $6$. Hence
$f$ fails $(T1)$, and its nonvanishing at order $6$ also gives a
failure of $(T2)$. Without the hypothesis $(T1)$ on $f$,
Question~\ref{q:weak} would therefore have a trivial negative answer.
The constructions below will retain this hypothesis.
\end{remark}

We use the definition of an unsupported divisor
from~\cite[Definition~1.5]{KLMS2026a}.
\begin{definition}\label{def:unsupported}
Let $A$ be a finite nonempty set and let
$N=\prod_{i=1}^k p_i^{\alpha_i}$, with $k\ge2$ and distinct primes $p_i$.
Then $\Ph_N$ is an \textbf{unsupported
divisor} of $A$ if $\Ph_N\mid A$, every $p_i$ divides $|A|$, and
$\Ph_{p_i^{\alpha_i}}\nmid A$ for every $i$.
The cardinality condition is part of this definition. For the
construction below, it will be enough to have the mixed zero and the
absence of its constituent prime-power zeros.
\end{definition}

\subsection{Normalized autocorrelations and positive definiteness}
The construction will start from the normalized autocorrelation of a
set complement. We recall why this gives a weak complement and how its
Fourier transform relates to that of the set.
\begin{definition}
A real-valued function $f:\zz_M\to\rr$ is \textbf{positive definite}
(PD) if $\widehat f(k)\ge0$ for every $k\in\zz_M$. A PD weak complement
satisfies this additional condition as well as \eqref{eq:weakintro}.
This is the finite-group notion of weak PD tiling
in~\cite[Definition~1.1]{KMMS2024}.
\end{definition}

\begin{samepage}
\begin{lemma}[G. Kiss, D. Matolcsi, M. Matolcsi, and G. Somlai, {\cite[Lemma~1.2]{KMMS2024}}]\label{lem:autocorrelation}
If $A\oplus B=\zz_M$, the autocorrelation function
\[
 g_B=\frac1{|B|}\one_B*\one_{-B}
\]
is a normalized PD weak complement of $A$ with mass $|B|$. Moreover,
\[
 \widehat{\one_B}(k)=0\quad\Longleftrightarrow\quad\widehat g_B(k)=0,
 \qquad k\in\zz_M.
\]
\end{lemma}
\end{samepage}

As the proof is direct and short, we include it for completeness.

\begin{proof}
Each coefficient is a difference count divided by the normalizing
factor, so it is nonnegative. Exactly $|B|$ ordered pairs of points of
$B$ have difference zero, giving the normalization. To check the
weak-tiling identity, associativity gives
$\one_A*g_B=|B|^{-1}\one_{\zz_M}*\one_{-B}=\one_{\zz_M}$.
Finally, the identity $\widehat g_B=|\widehat{\one_B}|^2/|B|$ proves
positive definiteness, preserves all Fourier zeros, and gives the mass.
\end{proof}

A real PD function is even: its Fourier transform is real, so Fourier
inversion and conjugation give $f(-x)=f(x)$. Later we will choose
perturbations that are not invariant under negation, and this observation
will show that the resulting weak complements are not PD.

\section{A pure-order cuboid filter}\label{sec:filter}
We use the cuboid formalism introduced in~\cite[Section~5]{LL2022}
and reviewed in~\cite[Section~2.3]{KLMS2026a}.
The construction below develops the cuboid perturbation proposed
in~\cite[Section~10.2]{KLMS2026a}: we construct a signed mask
that changes the Fourier values at a chosen order while leaving every
other Fourier value at the ambient period unchanged. Let $1<N\mid M$
and write
\[
 N=\prod_{i=1}^k p_i^{\alpha_i},\qquad \alpha_i\ge1,
\]
where the $p_i$ are distinct primes. Define the signed cuboid on $\zz_N$ by its
canonical mask
\[
 \Delta_N(X)=\prod_{i=1}^k(1-X^{N/p_i})\mod (X^N-1),
 \qquad \deg\Delta_N<N.
\]
Its $2^k$ vertices carry alternating signs. The \textbf{pure-order
cuboid filter} on $\zz_M$ is the periodic extension
\begin{equation}\label{eq:filter}
 R_{N,M}(X)=\Delta_N(X)\sum_{j=0}^{M/N-1}X^{Nj}.
\end{equation}
This polynomial has degree less than $M$; its coefficient function is
$R_{N,M}(x)=\Delta_N(x\bmod N)$.

\begin{samepage}
\begin{lemma}\label{lem:purefilter}
The filter $R=R_{N,M}$ has the following properties:
\begin{enumerate}[(1)]
\item $R(x)\in\{-1,0,1\}$, $R(0)=1$, and $\mass(R)=0$;
\item it has exactly $2^{k-1}M/N$ positive and the same number of negative
coefficients;
\item its Fourier support consists exactly of the characters of order $N$;
\item if $\Ph_N\mid A(X)$, then $\one_A*R=0$.
\end{enumerate}
\end{lemma}
\end{samepage}
\begin{proof}
To verify the coefficient claims, first consider the subset sums
$\sum_i\epsilon_iN/p_i$, $\epsilon_i\in\{0,1\}$. These are distinct
modulo $N$: reducing an equality modulo $p_i^{\alpha_i}$ eliminates
all terms except the $i$th, and the $p_i$-adic valuation of $N/p_i$
is $\alpha_i-1$, so the equality determines $\epsilon_i$. It follows
that the canonical mask polynomial $\Delta_N$ has $2^k$ nonzero coefficients,
each equal to $1$ or $-1$, with equally many of each sign and with
coefficient $1$ at the origin.

Repeating the cuboid does not change these coefficients. Indeed, the
shifted copies $X^{Nj}\Delta_N(X)$ in \eqref{eq:filter} occupy the
disjoint exponent intervals $[Nj,N(j+1)-1]$, so their supports do not
overlap. Each coefficient therefore occurs in a separate block,
proving (1) and (2).

For the Fourier support, let $\zeta$ have order $d\mid M$. The geometric sum
$\sum_j\zeta^{Nj}$ vanishes unless $d\mid N$, in which case it equals
$M/N$. If $d\mid N$ is proper, some $i$ satisfies $d\mid N/p_i$,
so the factor $1-\zeta^{N/p_i}$ vanishes. At order $N$, none of these
factors vanishes. This proves (3), including vanishing at the trivial
character. If $\Ph_N\mid A$, the Fourier transform of $\one_A$ vanishes
on the support of $\widehat R$, so the convolution theorem gives (4).
\end{proof}

\begin{samepage}
\begin{corollary}\label{cor:filteruse}
Let $Q\in\qq[X]$ and $t\in\qq\setminus\{0\}$. Adding
$tR_{N,M}$ leaves the values of $Q$ unchanged at every $M$th root
whose order differs from $N$. If $\Ph_N\mid Q$, then
$\Ph_N\nmid Q+tR_{N,M}$.

When $N$ is mixed, this change preserves every prime-power Fourier zero
and nonzero at the ambient period $M$.
\end{corollary}
\end{samepage}
\begin{proof}
By Lemma~\ref{lem:purefilter}(3), the filter vanishes at every $M$th
root except those of order $N$, where it is nonzero. The two assertions
follow by evaluating the perturbed polynomial at these roots.
\end{proof}

This gives a way to remove a mixed cyclotomic zero while retaining the
data entering the period-relative condition $(T1)$. To use it for weak
complements, we must also preserve nonnegativity and normalization.
Thus the negative filter sites need enough positive weight, and the
perturbation must leave the origin unchanged. We arrange these two
requirements by placing translated complements in separate layers.

\section{A cyclic construction without new primes}\label{sec:cycliclift}
We now use the filter to construct a normalized nonnegative weak
complement satisfying $(T1)$ and failing $(T2)$, starting from a mixed
zero of a set tile whose constituent prime-power zeros are absent.
The construction allows a dilation by a power of a prime already
dividing the mixed order, so no new prime is introduced into the period
or the coefficient denominators.

The use of the residue layers $j+q\zz_M$ of $\zz_{qM}$ is what permits
$q$ and $M$ to share prime factors. For comparison, the construction in
Appendix~\ref{sec:lift} uses a direct product $\zz_M\times\zz_r$,
with a new prime $r$ coprime to $M$.

\subsection{Removing off-period factors using an old prime}
Before constructing the weak complement, we show how to remove off-period
cyclotomic factors without changing the cyclic coefficient function.
We begin with the following coefficient reweighting result over $\qq[X]$.

\Needspace{12\baselineskip}
\begin{samepage}
\begin{lemma}\label{lem:oddtransfer}
Let $F\in\qq[X]$ have nonnegative coefficients and degree less than $D$.
Suppose that $mF\in\zz[X]$ for a positive integer $m$, and choose an integer
$x$ with $0<x<D$ and $[X^x]F>0$. For any odd prime $\ell$, the polynomial
\[
 P(X)=F(X)+\frac{X^{x+D}-X^x}{\ell m}
\]
has nonnegative coefficients, satisfies $P(0)=F(0)$ and $P(1)=F(1)$,
and represents the same cyclic coefficient function on $\zz_D$.
Every cyclotomic divisor of $P$ has order dividing $D$.
\end{lemma}
\end{samepage}
\begin{proof}
The only coefficient being decreased is the one at $x$. Since
$[X^x](mF)$ is a positive integer, we have $[X^x]F\ge1/m$, so
subtracting $1/(\ell m)$ preserves nonnegativity. The new coefficient
at $x+D$ is positive, and the two exponents agree modulo $D$.
Consequently the cyclic function is unchanged. As $x$ is nonzero,
the transfer also preserves the constant coefficient, and its two
terms cancel at one.

To locate the cyclotomic divisors, we work modulo the chosen odd prime.
The reduction identities used below also appear in
T.~Fallon, G.~Kiss, A.~Mayeli, and G.~Somlai~\cite[Lemma~2.3]{FKMS2026}
and R.~Malikiosis~\cite[proof of Theorem~6.1]{M2022}; we
include the details to keep track of root multiplicities.
Assume that $\Ph_n\mid P$. Since $\Ph_n$ is a monic polynomial, Gauss's lemma implies that $\Ph_n$ divides
the integer polynomial
\[
 \ell mP(X)=\ell mF(X)+X^x(X^D-1)
\]
in $\zz[X]$. For an integer polynomial $Q$, write $\overline Q$
for its coefficientwise reduction modulo $\ell$. Reducing the
divisibility gives
\[
 \overline{\Ph_n}(X)\mid X^x(X^D-1)
 \qquad\text{in }\mathbb F_\ell[X].
\]
Since $\overline{\Ph_n}(0)\ne0$, the polynomials $X^x$ and
$\overline{\Ph_n}$ are coprime, so
$\overline{\Ph_n}\mid X^D-1$ in $\mathbb F_\ell[X]$.
Write
\[
 D=\ell^bD_0,\qquad n=\ell^a n_0,\qquad
 (D_0,\ell)=(n_0,\ell)=1.
\]
This separates the order of a root from its multiplicity in the
chosen characteristic. Indeed, in characteristic $\ell$ we have
$X^D-1=(X^{D_0}-1)^{\ell^b}$, whose nonzero roots are the $D_0$th
roots of unity, each with multiplicity $\ell^b$.

Since $\ell\nmid n_0$, the polynomial $X^{n_0}-1$ is separable (i.e., it has $n_0$ distinct roots in an algebraic closure).
Its factorization into reduced cyclotomic polynomials shows that
$\overline{\Ph_{n_0}}$ has precisely the primitive $n_0$th roots,
each with multiplicity one. For $a\ge1$, the cyclotomic identity
\[
 \Ph_{\ell^a n_0}(X)\Ph_{n_0}(X^{\ell^{a-1}})
   =\Ph_{n_0}(X^{\ell^a})
\]
and the Frobenius identity
$Q(X^{\ell^j})=Q(X)^{\ell^j}$ for $Q\in\mathbb F_\ell[X]$ give
\[
 \overline{\Ph_{\ell^a n_0}}\,
 \overline{\Ph_{n_0}}^{\,\ell^{a-1}}
 =\overline{\Ph_{n_0}}^{\,\ell^a}.
\]
Cancellation in $\mathbb F_\ell[X]$ therefore yields
\[
 \overline{\Ph_{\ell^a n_0}}
 =\overline{\Ph_{n_0}}^{\,\ell^{a-1}(\ell-1)}.
\]
The reduced cyclotomic polynomials also have nonzero roots. Comparing
their orders with the roots allowed by the divisibility above forces
$n_0\mid D_0$. If $a\ge1$, comparing their multiplicities gives
\[
 \ell^{a-1}(\ell-1)\le\ell^b.
\]
Because $\ell$ is odd, $a\ge b+1$ would make the left-hand side
strictly larger than $\ell^b$. Hence $a\le b$, so $n\mid D$.
\end{proof}

The oddness of $\ell$ is needed only for the last multiplicity comparison.

\Needspace{16\baselineskip}
\subsection{The prime-preserving construction}
We now give our main result. Given a set tile $A$ of period $M$ that
admits an unsupported divisor $\Ph_N$ with $N\mid M$, one can construct a polynomial
$P\in\qq[X]$ with nonnegative coefficients such that
\[
 A(X^q)P(X)\equiv U_{qM}(X)\mod (X^{qM}-1),
\]
where $q=1$ or $q=p^h$ for a prime $p\mid N$ and an integer $h\ge1$.
The polynomial $P$ satisfies the Coven--Meyerowitz condition $(T1)$
but fails $(T2)$.
\begin{samepage}
\begin{theorem}\label{thm:primepreserving}
Suppose that $A\oplus B=\zz_M$ and that
\[
 N=\prod_{i=1}^k p_i^{\alpha_i}\mid M,\qquad k\ge2,
 \qquad \Ph_N\mid A,\qquad
 \Ph_{p_i^{\alpha_i}}\nmid A\quad(1\le i\le k).
\]
Then there are an integer $q\ge1$, which is either $1$ or a power of a
prime dividing $N$, and a polynomial $P\in\qq[X]$ such that the following hold.
\begin{enumerate}
\item
The coefficients of $P$ are nonnegative, $P(0)=1$, and
\[
 A(X^q)P(X)\equiv U_{qM}(X)\mod (X^{qM}-1).
\]
\item The polynomial $P$ satisfies the Coven--Meyerowitz condition $(T1)$ but fails $(T2)$.
\item If the coefficients $[X^x]P$ are written in lowest terms as $n_x/d_x \in \mathbb{Q}$, with $d_x\ge1$, then every prime divisor of $d_x$ divides $M$.
\item The dilated set $qA\subset\zz_{qM}$ is a set tile, has cardinality $|qA|=|A|$, and satisfies ordinary $(T2)$ if and only if $A$ does.
\end{enumerate}
\end{theorem}
\end{samepage}
\begin{proof}
Since $A \oplus B = \mathbb{Z}_M$, by Lemma~\ref{lem:allocation}, every $\Ph_{p_i^{\alpha_i}}$ divides $B$.
We prove (1)--(3) by considering whether the complement vanishes at
the same mixed order as the tile.

\smallskip\noindent
\emph{Case 1: $\Ph_N\nmid B$.}
Here the complement already has the required prime-power zeros and
lacks their mixed zero. We need only arrange that its indicator is
normalized. Choose $b_0\in B$ and replace $B$ by the translate
$B'=(B-b_0)\bmod M$, so that $0\in B'$ and
$A\oplus B'=\zz_M$. For every $d\mid M$, translation multiplies the
mask polynomial value at a primitive $d$th root by a nonzero monomial value.
It therefore preserves vanishing at every such order. Both set
complements have all their ordinary prime-power factors at orders
dividing $M$, by Lemma~\ref{lem:allocation}; hence their ordinary
prime-power signatures agree. Every mixed order required by $(T2)$
also divides $M$, so translation preserves $(T2)$ and its failure.

Set $q=1$ and let $P$ be the canonical mask polynomial of $B'$.
Then $P$ has nonnegative integer coefficients,
$P(0)=1$, and $A(X)P(X)\equiv U_M(X)\mod (X^M-1)$.
It satisfies ordinary $(T1)$ by allocation. All the constituent
prime-power factors of $N$ divide $P$, whereas $\Ph_N\nmid P$,
so $(T2)$ fails at $N$.

\smallskip\noindent
\emph{Case 2: $\Ph_N\mid B$.}
In this case both mask polynomials vanish at the mixed order. We will build a
layered complement whose autocorrelation has enough positive weight
at every negative filter site. Let $R=R_{N,M}$ be the filter of
Lemma~\ref{lem:purefilter}; its negative support has cardinality
\[
 s=2^{k-1}M/N.
\]
List these positions as $x_1,\ldots,x_s\in\{0,\ldots,M-1\}$.
We allow a consecutive pair of layers for each of these positions by
choosing a prime $p\mid N$ and a power $q=p^h$ such that
$q\ge\max(3,s+1)$. Put $D=qM$ and $K=q|B|$. Choose shifts in $\zz_M$ by
\[
 t_0=0,\qquad t_j-t_{j-1}=x_j\quad(1\le j\le s),\qquad
 t_j-t_{j-1}=0\quad(s<j<q).
\]
The layered set
\begin{equation}\label{eq:newlayers}
 C=\bigcup_{j=0}^{q-1}\{j+q(b+t_j):b\in B\}\subset\zz_D
\end{equation}
has $K$ elements and complements $qA$. Indeed, its layers occupy
distinct residue classes modulo $q$. Subtracting $j$ and dividing
by $q$ identifies the tiling in layer $j$ with
$A\oplus(B+t_j)=\zz_M$. Equivalently,
\[
 C(X)\equiv B(X^q)\sum_{j=0}^{q-1}X^{j+qt_j}\mod (X^D-1).
\]
Multiplication by $A(X^q)$ gives $U_M(X^q)U_q(X)=U_D(X)$,
because $U_M(X^q)$ is invariant under multiplication by $X^{qt_j}$
in the quotient ring $\qq[X]/(X^D-1)$.

To keep the origin fixed, we place the filter in a nonzero residue
layer of the autocorrelation. Define functions on $\zz_D$ by
\begin{equation}\label{eq:newperturbation}
 g=\frac1K\one_C*\one_{-C},\qquad
 H(X)=X R(X^q),\qquad f=g+\frac1K H.
\end{equation}
The polynomial $H$ has degree less than $D$, and its coefficient function
is supported in the residue class $1\bmod q$. Explicitly,
$H(1+qx)=R(x)$ for $0\le x\le M-1$. Since $\one_A*R=0$, substitution gives
\[
 A(X^q)H(X)=X A(X^q)R(X^q)=X(AR)(X^q)\equiv0\mod (X^D-1).
\]
Also $H(0)=0$ and $\mass(H)=0$.

We now address the nonnegativity of the weight function $f$. For each negative position $x_j$ and every $b\in B$, the two points
\[
 j+q(b+t_j),\qquad (j-1)+q(b+t_{j-1})
\]
of $C$ have difference $1+qx_j$ modulo $D$. These are $|B|$ distinct
ordered pairs, so
\[
 g(1+qx_j)\ge\frac{|B|}{K},\qquad
 f(1+qx_j)\ge\frac{|B|-1}{K}\ge0.
\]
These are all the negative sites of $H$, and consequently $f\ge0$.
Lemma~\ref{lem:autocorrelation} gives $g(0)=1$,
$\one_{qA}*g=\one_{\zz_D}$, and $\mass(g)=K$. Since
$H(0)=\mass(H)=0$ and $\one_{qA}*H=0$, the same identities hold for $f$:
\[
 f(0)=1,\qquad \one_{qA}*f=\one_{\zz_D},\qquad \mass(f)=K.
\]

We next check that the perturbation preserves the prime-power zeros
while removing the mixed zero. For a root of unity $\zeta$ of order
dividing $D$, the value $H(\zeta)$
is nonzero exactly when $\zeta^q$ has order $N$. If $\zeta$ has order
$d$, then $\zeta^q$ has order $d/\gcd(d,q)$. Since $q$ is a power
of a prime dividing $N$, this order equals $N$ exactly when $d=qN$.
Thus $H$ has pure Fourier order $qN$. As $qN$ has at least two
distinct prime factors, $H$ vanishes at every prime-power character.
The prime-power Fourier zeros of $f$ are therefore exactly those of
$g$, hence of $C$. Lemma~\ref{lem:allocation} and $|C|=K$ prove
period-relative $(T1)$.

At a primitive $qN$th root $\zeta$, the layered factorization of $C$
gives $C(\zeta)=0$, since $B(\zeta^q)=0$. The autocorrelation $g$
therefore also vanishes at this character. In contrast,
$H(\zeta)=\zeta R(\zeta^q)\ne0$, so $f$ does not vanish at order $qN$.
To obtain the $(T2)$ failure, we also need its constituent prime-power
zeros. The constituent prime powers of $qN$ are
$p_i^{\alpha_i+v_{p_i}(q)}$. Evaluating $A(X^q)$ at each such order
reduces to a nonzero evaluation of $A$ at order $p_i^{\alpha_i}$.
Allocation of prime-power divisors puts all these prime-power factors in $C$, and hence in $f$. This then proves the $(T2)$ failure.

It remains to pass from the cyclic function to an ordinary polynomial
with the required conditions. Let $F$ be the canonical mask polynomial of $f$.
We have $KF\in\zz[X]$ and
$[X]F\ge1/K$, because $[X]H=[X^0]R=1$. Since $k\ge2$, some odd
prime $\ell$ divides $N$. Lemma~\ref{lem:oddtransfer}, with $m=K$
and $x=1$, gives the ordinary representative
\begin{equation}\label{eq:newrepresentative}
 P(X)=F(X)+\frac{X^{D+1}-X}{\ell K}.
\end{equation}
By Lemma~\ref{lem:oddtransfer}, $P$ has nonnegative coefficients, satisfies $P(0) = 1$, and is congruent to $F$ modulo $X^{D} - 1$, proving (1).
It has no off-period cyclotomic divisors. Thus period-relative $(T1)$
becomes ordinary $(T1)$, while the cyclic normalization and $(T2)$
failure remain unchanged. This proves (2) in Case 2. Every prime in the denominator
$\ell K=\ell q|B|$ divides the original period $M$, verifying (3).

To prove (4), observe that for every $q\ge1$, evaluation at prime-power roots gives
\begin{equation}\label{eq:dilatedsignature}
 S_{A(X^q)}^*=
 \{\ell^{a+v_\ell(q)}:\ell^a\in S_A^*\}.
\end{equation}
At order $\ell^b$ with $b\le v_\ell(q)$, the value is $|A|$, not zero;
at larger orders, taking the $q$th power reduces the root order by
the factor $\ell^{v_\ell(q)}$. For products of powers of distinct
primes, the same calculation applies to all constituent exponents
simultaneously. The mixed orders required by $(T2)$ for $A(X^q)$
therefore correspond exactly to those required for $A$, and evaluation
at their primitive roots proves $(T2)$ in both directions.
The set $qA$ has cardinality $|A|$ and is a tile: the complement is
$B'$ in Case~1 and $C$ in Case~2.
\end{proof}

\begin{remark}
In any nontrivial-dilation branch where $q = p^h$, $g$ is even whereas the cyclic coefficients of $H$ at $1$ and $-1$ are respectively $1$ and $0$,
since $q>2$. Thus $f(1)-f(-1)=1/K$, and $f$ is not positive definite.
\end{remark}

The next remark explains how a set complement failing $(T2)$ would
yield a PD weak complement with the same failure. If the set factor
also satisfies $(T2)$, this would give a counterexample to the
positive-definite analogue of Question~\ref{q:weak}, and hence to the
extension proposed in~\cite[Problem~1.2]{KLMS2025}. No such set
complement is constructed here.

\begin{remark}\label{rem:newnonpd}
Case~1 begins with a set complement $B$ that already fails $(T2)$.
If such a complement exists, its normalized autocorrelation is a PD
weak complement with the same Fourier zeros and the same failure.
Its mass is $|B|>1$ and its value at zero is $1$, so it has a positive
coefficient away from zero. Applying Lemma~\ref{lem:oddtransfer}
with $m=|B|$ gives an ordinary $(T1)$ representative without
changing that cyclic PD function. Thus Case~1 yields a PD weak complement
failing $(T2)$ only if one already has a set complement failing $(T2)$.
If $A$ also satisfies $(T2)$, this would answer the positive-definite
analogue of Question~\ref{q:weak} negatively.
\end{remark}
\section{A three-prime counterexample of period 180}\label{sec:example180}

We now prove Theorem~\ref{thm:introexample} by giving the weak complement
explicitly. We start from the five-point set $\{0,2,3,4,6\}$, whose
mask polynomial is $\Ph_5\Ph_6$, and follow the construction of the preceding
section with mixed order $N=6$. In this example, six translated layers
suffice to keep the perturbed complement nonnegative. This gives the composite
dilation $q=6$, although Theorem~\ref{thm:primepreserving} chooses a
prime power and uses the larger bound $q\ge s+1=11$. The resulting
period is $180$, with only three distinct prime divisors. The set
realizes~\cite[Example~8.1]{KLMS2026a} with primes $2,3,5$, as noted in
the introduction; we verify below all of its properties that we use.

For a positive integer $s$, write $U_s(Y)=1+Y+\cdots+Y^{s-1}$. Define
\begin{equation}\label{eq:180sets}
 A_* =\{0,12,18,24,36\}\subset\zz_{180},\qquad
 T=\{0,1,5,8,10,21\}\subset\zz_{30}.
\end{equation}
Let $Q\in\zz[X]$ be the canonical cyclic autocorrelation mask polynomial of $T$:
\begin{equation}\label{eq:180Q}
 Q(X)=\sum_{t,t'\in T}X^{(t-t')\bmod30},\qquad
 0\le (t-t')\bmod30<30.
\end{equation}
Using $Q$, define the following ordinary polynomials:
\begin{align}
 H(X)&=X(1-X^{12})(1-X^{18})U_5(X^{36}),\label{eq:180H}\\
 W(X)&=Q(X)U_6(X^{30})+H(X),\qquad F(X)=\frac16W(X).
 \label{eq:180F}
\end{align}
The polynomials $H$, $W$, and $F$ all have degree less than $180$.

\Needspace{16\baselineskip}
\begin{samepage}
\begin{theorem}\label{thm:180}
The set $A_*$ defined in \eqref{eq:180sets} tiles $\zz_{180}$ and
satisfies ordinary $(T1)$ and $(T2)$.

The polynomial $F$ defined in \eqref{eq:180F} has degree $179$,
nonnegative rational coefficients, and $F(0)=1$, and it satisfies
\[
 A_*(X)F(X)\equiv U_{180}(X)\mod (X^{180}-1).
\]
Moreover,
\[
 F(1)=36,\qquad S_F^*=\{2,3,4,9\},\qquad \Ph_{36}\nmid F.
\]
Thus $F$ satisfies ordinary $(T1)$ and fails $(T2)$.
Its cyclic coefficient function is not positive definite.
\end{theorem}
\end{samepage}

\begin{proof}
\emph{The base tiling and the lifted set complement.}
We first check that $A_*$ tiles. An explicit subgroup complement is
\[
 B_*=5\zz_{180}=\{0,5,10,\ldots,175\},\qquad |B_*|=36.
\]
Indeed, the elements of $A_*$ represent the five residue classes modulo
$5$ exactly once, and hence $A_*\oplus B_*=\zz_{180}$. To construct
the required weak complement, we use the autocorrelation of a different
set complement $C$, whose layers we choose as follows.

In $\zz_{30}$ take
\[
 A_0=\{0,2,3,4,6\},\qquad B_0=\{0,5,10,15,20,25\}.
\]
The set $A_0$ represents each residue class modulo $5$ exactly once, so
$A_0\oplus B_0=\zz_{30}$. Direct multiplication gives
\[
 A_0(X)=1+X^2+X^3+X^4+X^6=\Ph_5(X)\Ph_6(X).
\]
Choose the six shifts
\[
 (t_0,t_1,t_2,t_3,t_4,t_5)=(0,0,1,3,6,10)
\]
and define, with all residues reduced modulo $180$,
\begin{equation}\label{eq:180C}
 C=\bigcup_{j=0}^5\{j+6(b+t_j):b\in B_0\}.
\end{equation}
To see that $C$ complements $A_*$, observe that its six layers lie in
distinct residue classes modulo $6$. In the class indexed by $j$,
subtracting $j$ and dividing by $6$ reduces the required tiling identity
to $A_0\oplus(B_0+t_j)=\zz_{30}$. Thus $A_*\oplus C=\zz_{180}$,
with $|C|=36$.

Modulo $30$, the residues $j+6t_j$ are $0,1,8,21,10,5$.
Consequently
\[
 C=T+\{0,30,60,90,120,150\},\qquad
 C(X)=T(X)U_6(X^{30}).
\]
Writing $q_n=[X^n]Q$ for $0\le n<30$, we obtain
\begin{equation}\label{eq:180autocorrelation}
 (\one_C*\one_{-C})(n)=6q_{n\bmod30},\qquad n\in\zz_{180}.
\end{equation}
For each ordered pair of representatives from $T$, there are six
ordered pairs from $C$ with any prescribed difference in the
corresponding coset modulo $30$. This proves
\eqref{eq:180autocorrelation} and shows that the first term of $F$ in
\eqref{eq:180F} is the canonical mask polynomial of the normalized autocorrelation
$g=(\one_C*\one_{-C})/36$.

\emph{The signed perturbation and nonnegativity.}
On $\zz_{30}$ let
\[
 R(X)=(1-X^2)(1-X^3)U_5(X^6).
\]
The cyclic coefficient function of $R$ repeats the vector
$(1,0,-1,-1,0,1)$ five times, and its Fourier transform is supported
exactly at characters of order $6$. Indeed, $U_5(X^6)$ vanishes at
the $30$th roots whose orders do not divide $6$, while the remaining
two factors vanish at orders $1$, $2$, and $3$, but not at order $6$.
Since $\Ph_6\mid A_0$, it follows that
\[
 A_0(X)R(X)\equiv0\mod (X^{30}-1).
\]
The polynomial \eqref{eq:180H} is exactly $XR(X^6)$. Hence
\begin{equation}\label{eq:180annihilation}
 A_*(X)H(X)=X A_0(X^6)R(X^6)
      \equiv0\mod (X^{180}-1).
\end{equation}
The exponents occurring in $H$ lie in the residue class $1\bmod 6$,
and $[X^{1+6x}]H=[X^x]R\in\{-1,0,1\}$ for $0\le x<30$.
To check that the negative coefficients can be absorbed by the
first term of $W$, we count the corresponding differences in $T$:
\[
 \begin{array}{c|ccccc}
 x\bmod5&0&1&2&3&4\\ \hline
 q_{(1+6x)\bmod30}&1&1&1&1&2.
 \end{array}
\]
The entries in the second row are all at least one, so every negative
coefficient of $H$ is absorbed and $W$ has nonnegative coefficients.
Since the constant coefficient of $H$ is zero and $q_0=|T|=6$, we also
have $F(0)=1$. Moreover, $Q(1)=|T|^2=36$ and $H(1)=0$, from which
$W(1)=216$ and $F(1)=36$ follow. Finally, the set tiling gives
$\one_{A_*}*g=\one_{\zz_{180}}$, and
\eqref{eq:180annihilation} shows that the perturbation preserves this
identity. This proves the required weak tiling.

\emph{The prime-power zeros and ordinary $(T1)$.}
Multiplication by $6$ permutes the residue classes modulo $5$, so
$\Ph_5\mid A_*$. Since $A_*(1)=5$, this is its only prime-power
cyclotomic factor: a second such factor would force an additional prime
factor in the integer value $A_*(1)$. This proves ordinary $(T1)$ for
$A_*$, and $(T2)$ has no nontrivial mixed requirement.
The weak-tiling identity now forces
\[
 \Ph_2\Ph_4\Ph_3\Ph_9\mid W.
\]
It remains to exclude any further prime-power divisors, including
those whose orders do not divide $180$. Since $W$ has integer
coefficients and $W(1)=216=2^3 3^3$, a divisor $\Ph_{p^a}$ must have
$p\in\{2,3\}$. Furthermore, $q_{29}=1$ and $\deg H=175$, so that
$\deg W=179$. Comparing this with the degree of $\Ph_{p^a}$ leaves
only the orders
\[
 2,4,8,16,32,64,128,256,\quad \textrm{and} \quad  3,9,27,81,243.
\]
The constant coefficients of the monic remainders at the nine unwanted
orders are as follows:
\begin{equation}\label{eq:180exclusions}
\begin{array}{c|rrrrrrrrr}
 d&8&16&32&64&128&256&27&81&243\\ \hline
 [X^0](W\bmod\Ph_d)&6&2&6&5&6&5&6&8&6.
\end{array}
\end{equation}
The table can be checked by hand as follows.
Put $w_n=[X^n]W$ for $0\le n<180$ and $w_n=0$ otherwise. For $d=p^a$,
the entry in the second row is
\[
 \sum_{k\ge0}\bigl(w_{kd}-w_{(p-1)p^{a-1}+kd}\bigr).
\]
This follows by first reducing modulo $X^d-1$ and then using
$\Ph_{p^a}=1+X^{p^{a-1}}+\cdots+X^{(p-1)p^{a-1}}$.
Each sum is finite, and every entry of
\eqref{eq:180exclusions} is nonzero. Consequently
$S_F^*=\{2,3,4,9\}$, and the values of the corresponding cyclotomic
polynomials at one have product $2\cdot2\cdot3\cdot3=36=F(1)$.
This proves ordinary $(T1)$.

\emph{Failure of $(T2)$ and positive definiteness.}
At a primitive $36$th root, the geometric sum $U_6(X^{30})$ vanishes.
Since
\[
 \Ph_{36}(X)=X^{12}-X^6+1,
 \qquad X^{18}\equiv-1\mod \Ph_{36},
\]
formula \eqref{eq:180H} gives
\[
 W\equiv H\equiv 5X(1-X^{12})(1-X^{18})
       \equiv20X-10X^7\mod \Ph_{36}.
\]
The remainder is nonzero and has degree less than
$\deg\Ph_{36}=12$, proving that $\Ph_{36}\nmid F$. Since both
$\Ph_4$ and $\Ph_9$ divide $F$, this gives the failure of $(T2)$.
Finally, $w_1=q_1+R(0)=2$, whereas $w_{179}=q_{29}=1$. Thus the cyclic
coefficient function of $F$ is not even, and hence is not positive
definite.
\end{proof}

\begin{remark}
The weak complement $F$ belongs to the dilated set $A_*=6A_0$, not to
$A_0$ itself. The next section shows that this is unavoidable in general:
some tiles with an unsupported divisor have no such weak complement at
any period.
\end{remark}

\section{What period enlargement does and does not preserve}\label{sec:obstructions}
The construction of Theorem~\ref{thm:primepreserving} replaces $A(X)$ by
$A(X^q)$ before constructing a weak complement. We now examine the role
of this change in the set factor. First, we give a tile for which
enlarging the ambient period, while keeping its mask polynomial fixed, can never
produce a normalized nonnegative $(T1)$ weak complement failing $(T2)$.
We then show that dilation cannot create a set complement failing
$(T2)$, and that the rational construction requires a period with at
least three distinct prime divisors.

\subsection{Allocation for a rational complement satisfying \texorpdfstring{$(T1)$}{(T1)}}
We begin with the allocation of prime-power factors. For a rational
complement, integrality no longer gives the counting argument used in
Lemma~\ref{lem:allocation}. The condition $(T1)$ supplies the missing
constraint through the value of the complement at one.

\begin{samepage}
\begin{lemma}\label{lem:weakallocation}
Let $A\in\zz[X]$ satisfy $a=A(1)>0$, and suppose that $P\in\qq[X]$
satisfies ordinary $(T1)$ and
\[
 A(X)P(X)\equiv U_M(X)\mod (X^M-1).
\]
Then $A$ satisfies ordinary $(T1)$, and $S_A^*$ and $S_P^*$ form a
disjoint partition of the prime-power orders dividing $M$.
No positivity assumption is needed.
\end{lemma}
\end{samepage}
\begin{proof}
We compare the prime-power factors forced into $P$ by the cyclic
identity with those already present in $A$. Put
\[
 Q_A=\prod_{\substack{p^b\mid M\\\Ph_{p^b}\mid A}}p.
\]
The product of the distinct monic cyclotomic factors counted in $Q_A$
has an integral quotient in $A$. Evaluating at one gives $Q_A\mid a$,
and in particular $Q_A\le a$. On the other hand, each prime-power order
dividing $M$ that is absent from $A$ must occur in $P$. The values at one
of these forced factors multiply to $M/Q_A$. Since the cyclic identity
also determines $P(1)$, ordinary $(T1)$ for $P$ gives
\[
 \frac{M}{a}=P(1)
 =\prod_{s\in S_P^*}\Ph_s(1)\ge\frac{M}{Q_A}.
\]
It follows that $Q_A\ge a$, so both inequalities must be equalities.
The prime-power factors already counted in $A$ therefore exhaust its
value at one, leaving no further prime-power factor. Equality for $P$
likewise excludes both a factor shared with $A$ and a factor whose
order does not divide $M$. This proves the asserted partition and
$(T1)$ for $A$.
\end{proof}

\subsection{An unchanged tile obstructed at every period}
The fixed-tile question is motivated by the discussion
in~\cite[Section~10.2]{KLMS2026a}. To show why the set factor cannot in
general be left unmodified (and, in our specific context, undilated), we use the period-$180$ example of H.~Tan
and Y.~Zhang~\cite[Section~8]{TZ2026}.

This tile consists of six translates of the five-point set from
Section~\ref{sec:example180}. Its nonnegative normalized complements
satisfy a rigidity property that persists when the ambient period is
enlarged.

\begin{samepage}
\begin{proposition}\label{prop:fixedmask}
Let
\[
 E=\{0,2,3,4,6\},\qquad
 T=E+30\{0,1,2,3,4,5\}.
\]
The set $T$ tiles $\zz_{180}$, satisfies ordinary $(T1)$ and $(T2)$,
and has unsupported divisor $\Ph_6$. Nevertheless, for every integer
$L\ge2$, every polynomial $P\in\qq[X]$ with nonnegative coefficients
satisfying
\[
 P(0)=1,\qquad
 T(X)P(X)\equiv U_L(X)\mod (X^L-1),
\]
and ordinary $(T1)$ must also satisfy ordinary $(T2)$.
\end{proposition}
\end{samepage}
\begin{proof}
We first verify the tiling and cyclotomic data. The set $E$ represents
each residue class modulo five exactly once, so
$E\oplus B=\zz_{30}$ for $B=5\{0,1,2,3,4,5\}$. Consequently
$T\oplus B=\zz_{180}$. Direct multiplication gives
\begin{equation}\label{eq:fixedmaskfactorization}
 T(X)=\frac{X^{180}-1}{X^{30}-1}\Ph_5(X)\Ph_6(X),
 \qquad E(X)=\Ph_5(X)\Ph_6(X).
\end{equation}
In particular, $S_T^*=\{4,5,9\}$ and $|T|=30$. Every mixed product of
prime powers in this signature divides $180$ but not $30$, and hence
occurs among the factors in the quotient in
\eqref{eq:fixedmaskfactorization}. This proves $(T1)$ and $(T2)$.
The same factorization gives $\Ph_6\mid T$ and
$\Ph_2,\Ph_3\nmid T$. Since $2,3\mid |T|$, the divisor $\Ph_6$ is
unsupported.

The main point of the proof is that the set $E$ has only one
normalized nonnegative weak complement on $\zz_{30}$. To see this,
suppose that
$F\ge0$ on $\zz_{30}$, $F(0)=1$, and
$\one_E*F=\one_{\zz_{30}}$. Its mass is six. By
\eqref{eq:fixedmaskfactorization}, the nontrivial Fourier support of
$F$ can occur only at orders five and six. We may therefore separate
these two components in the CRT coordinates
$\zz_{30}\simeq\zz_5\times\zz_6$. Assigning the constant Fourier term
to the first component, Fourier inversion gives real functions $u$
and $v$ such that
\[
 F(x,y)=u(x)+v(y),\qquad
 \sum_{x\in\zz_5}u(x)=1,\qquad
 \sum_{y\in\zz_6}v(y)=0,
\]
where $v$ has Fourier support only at order six. Put
$c=-\min_y v(y)\ge0$. Nonnegativity of $F$ gives $u(x)\ge c$ for every
$x$, and normalization gives
\[
 v(0)=1-u(0)=\sum_{x\ne0}u(x)\ge4c.
\]
The order-six support imposes the additional relation
$v(0)+v(2)+v(4)=0$, from which $v(0)\le2c$. The two bounds on $v(0)$
are compatible only if $c=0$. Thus $v\ge0$, and its zero mass implies
$v=0$. We then have $u(0)=1$ and $u(x)=0$ for $x\ne0$, proving that
the only normalized nonnegative weak complement of $E$ on $\zz_{30}$ is
\begin{equation}\label{eq:foldrigidity}
 F=\one_{5\zz_{30}}.
\end{equation}

We next transfer this rigidity to complements of $T$. Suppose that
$h\ge0$ on $\zz_{180}$, $h(0)=1$, and
$\one_T*h=\one_{\zz_{180}}$, and fold $h$ modulo $30$ by setting
\[
 F(y)=\sum_{j=0}^5h(y+30j),\qquad y\in\zz_{30}.
\]
The six translates in the definition of $T$ give
$\one_E*F=\one_{\zz_{30}}$. Moreover, $F(0)\ge h(0)=1$, whereas the
convolution equation at zero, together with $0\in E$ and
nonnegativity, gives $F(0)\le1$. We may therefore apply
\eqref{eq:foldrigidity}. Since folding preserves Fourier values at
characters whose orders divide $30$, it follows that the Fourier
transform of $h$ vanishes at order six. This conclusion required no
$(T1)$ assumption on $h$.

It remains to show that a larger ambient period gives no additional
freedom relevant to $(T2)$. Let $P$ satisfy the stated hypotheses at
period $L$. By Lemma~\ref{lem:weakallocation},
\begin{equation}\label{eq:fixedallocation}
 180\mid L,\qquad
 S_P^*=\{p^a:p^a\mid L\}\setminus\{4,5,9\}.
\end{equation}
Write $f=w_P^L$. The factorization
\eqref{eq:fixedmaskfactorization} shows that every cyclotomic factor of
$T$ has order dividing $180$. The cyclic identity therefore confines
the Fourier support of $f$ to characters of these orders, together
with the trivial character, and hence $f(n+180)=f(n)$. We must also
check that this $180$-periodic pattern retains the normalization.
Nonnegative coefficients and $P(0)=1$ give
$f(0)\ge1$. The convolution equation at zero, using $0\in T$, gives the
reverse inequality. The $180$-periodic pattern $h(n)=f(n)$, $0\le n<180$,
therefore satisfies
\[
 h\ge0,\qquad h(0)=1,\qquad
 \one_T*h=\one_{\zz_{180}}.
\]
Applying the preceding argument to $h$ gives the required order-six
vanishing, and therefore $\Ph_6\mid P$.

We can now verify $(T2)$ using the complete list of cyclotomic factors
of $T$. Let $d>1$ be a product of prime-power orders of distinct primes
in $S_P^*$. By \eqref{eq:fixedallocation}, $d\mid L$. If $\Ph_d\nmid P$,
the cyclic identity forces $\Ph_d\mid T$. By
\eqref{eq:fixedmaskfactorization}, the possible orders are $5$, $6$,
and the divisors of $180$ that do not divide $30$. Order five is excluded
by \eqref{eq:fixedallocation}; every order in the last class has
constituent prime power four or nine, also excluded. At the remaining
order six, we have already proved $\Ph_6\mid P$. This contradiction
proves $(T2)$.
\end{proof}

The period reduction used in this proof reflects a general restriction
on a fixed mask. For a nonzero integer mask polynomial $A$, let $L_A$ be the
least common multiple of the orders of its cyclotomic factors, with
$L_A=1$ if there are none. At ambient period $L$, the Fourier support of
every cyclic complement of $A$ has orders dividing $\gcd(L,L_A)$.
The complement consequently has a period dividing this number.
Thus enlarging the period cannot introduce Fourier orders beyond
those already present among the cyclotomic divisors of $A$. The allocation and normalization arguments
above are additionally needed for the stated $(T1)$--$(T2)$ obstruction.

\subsection{Dilation and set complements}
For set complements, the residue layers of a dilation give a further
restriction: a failure of $(T2)$ in a complement of $qA$ must already
occur in one of its layers. The following argument uses only allocation
of prime-power factors and does not depend on the integer tiling
characterization claimed by C.~Li et al.~\cite{LWHWWZ2026}.

\begin{samepage}
\begin{proposition}\label{prop:dilationequivalence}
Let $A\subset\zz_M$ tile $\zz_M$ and let $q\ge1$. Then $A$ has a
set complement failing $(T2)$ if and only if
$qA\subset\zz_{qM}$ has such a complement.
\end{proposition}
\end{samepage}
\begin{proof}
A translation by an element of $qA$ preserves each residue class modulo
$q$. Thus every set complement $C$ of $qA$ decomposes uniquely into
residue layers
\[
 C=\bigcup_{j=0}^{q-1}(j+qB_j),\qquad
 A\oplus B_j=\zz_M\quad(0\le j<q).
\]
We first show that $(T2)$ for every layer implies $(T2)$ for $C$.
Suppose, accordingly, that all the $B_j$ satisfy $(T2)$. Their ordinary
prime-power signatures agree by Lemma~\ref{lem:allocation}. Let $d>1$
be a product of prime-power orders of distinct primes in $S_C^*$,
choose a primitive $d$th root $\zeta$, and put
$e=d/\gcd(d,q)$. Then $\zeta^q$ has order $e$ and
\[
 C(\zeta)=\sum_{j=0}^{q-1}\zeta^j B_j(\zeta^q).
\]
If $e=1$, then $d\mid q$, and the values of the layer mask polynomials in this
sum are simply their cardinalities. Each is equal to $M/|A|$, so the
sum becomes $(M/|A|)U_q(\zeta)=0$.

For $e>1$, we use the common prime-power signature of the layers.
A constituent $p^{\beta-v_p(q)}$ of $e$ comes from a
constituent $p^\beta$ of $d$ with $\beta>v_p(q)$. Since
$p^\beta\in S_C^*$, allocation places this order outside the signature
of $A(X^q)$. Formula~\eqref{eq:dilatedsignature} then shows that
$p^{\beta-v_p(q)}$ is absent from $S_A^*$ and belongs to every
$S_{B_j}^*$. By $(T2)$ for the layers, or by their prime-power vanishing
if $e$ is a prime power, $B_j(\zeta^q)=0$ for every $j$. Thus $C$
satisfies $(T2)$. Contrapositively, a complement of $qA$ that fails $(T2)$
has a layer that is already a complement of $A$ failing $(T2)$.

For the converse, we use identical layers. Suppose that $B$ complements
$A$ and fails $(T2)$ at order $e$. The set with mask
\[
 C(X)=U_q(X)B(X^q)
\]
complements $qA$. Put $d=e\prod_{p\mid e}p^{v_p(q)}$. The constituent
prime powers of $d$ belong to $S_C^*$ by allocation and
\eqref{eq:dilatedsignature}. At a primitive $d$th root $\zeta$, the
root $\zeta^q$ has order $e$, so $B(\zeta^q)\ne0$. Also
\[
 U_q(\zeta)=\frac{1-\zeta^q}{1-\zeta}\ne0,
\]
because $e>1$. Hence $C(\zeta)\ne0$, giving its $(T2)$ failure.
\end{proof}

The rational construction avoids this set-level obstruction by adding
a fractional perturbation to a normalized autocorrelation. Its weak
complement need not be an indicator in any residue layer.

\Needspace{10\baselineskip}\subsection{The two-prime boundary}
Finally, we discuss what happens when the tiling period $M$ has at most two distinct prime divisors. As the following simple result shows, no such nonnegative normalized weak tiling complement satisfying $(T1)$ but failing $(T2)$ can exist.

\begin{samepage}
\begin{proposition}\label{prop:twoprime}
Suppose that $M$ has at most two distinct prime divisors and
$A\oplus B=\zz_M$. If $P\in\qq[X]$ satisfies ordinary $(T1)$ and
\[
 A(X)P(X)\equiv U_M(X)\mod (X^M-1),
\]
then $P$ satisfies ordinary $(T2)$. Positivity and normalization of
$P$ are not required.
\end{proposition}
\end{samepage}
\begin{proof}
By Lemma~\ref{lem:weakallocation}, the prime-power orders dividing $M$
are partitioned between $A$ and $P$. If $M$ has only one prime divisor,
there is no mixed requirement. Otherwise, write its two prime divisors
as $p$ and $r$, and suppose that $P$ fails $(T2)$ at
$N=p^\alpha r^\beta\mid M$. Both constituent prime powers must then
be absent from $A$, whereas the cyclic identity forces $\Ph_N\mid A$.

Since $|A|\mid M$, the two-prime theorem of E.~Coven and A.~Meyerowitz~\cite{CM1999}
gives $(T2)$ for $A$. By G.~Kiss, I.~{\L}aba, C.~Marshall, and
G.~Somlai~\cite[Theorem~9.4]{KLMS2026a}, a
nonnegative integer multiset on a two-prime cyclic group satisfying
$(T2)$, with a mixed zero but neither constituent prime-power zero,
has mass strictly larger than the product of its prime-power
cyclotomic values at one. Applied to $A$, this contradicts $(T1)$.
This version of the cited theorem does not require the additional
cardinality condition in Definition~\ref{def:unsupported}.
\end{proof}

A period with at least three distinct prime divisors is therefore
necessary for a rational complement satisfying $(T1)$ to fail $(T2)$. Without
the $(T1)$ assumption this restriction disappears, as the example in
$\zz_6$ following Lemma~\ref{lem:allocation} shows.

\section{Positive definiteness and the present constructions}\label{sec:pdlimitations}
Every spectral set in a finite cyclic group admits a normalized
nonnegative weak complement that is positive
definite~\cite[Lemma~1.3]{KMMS2024}.
G.~Kiss, I.~Londner, M.~Matolcsi, and
G.~Somlai~\cite[Problems~1.1 and~1.2]{KLMS2025} ask whether every set
in a finite cyclic group admitting a PD weak complement is a tile,
and whether the Coven--Meyerowitz conditions extend to such tilings.
They answer the first question affirmatively for prime-power
periods~\cite[Proposition~4.1]{KLMS2025}. They also construct functional
PD tilings $f*g=\one_{\zz_M}$, with $M=p^4q^2$ for primes $p<q<p^2$,
in which both factors satisfy $(T1)$ and neither satisfies
$(T2)$~\cite[Theorem~1.6]{KLMS2025}. Neither factor is the indicator
function of a set. Related weak-tiling questions are considered
in~\cite{KF2026,KLMS2026b}. Geometric consequences of weak tiling in
Euclidean spaces are studied in~\cite{KLM2025}.

The explicit weak complements were shown to fail positive definiteness
by comparing their values at opposite points. There is also a Fourier
reason why the perturbation cannot preserve this property. We start
with an autocorrelation, change only Fourier values that are initially
zero, and leave the value at the origin unchanged. The following
proposition shows that these requirements exclude every nonzero
perturbation on that Fourier support.

\begin{samepage}
\begin{proposition}\label{prop:pdperturbation}
Let $g,H:\zz_D\to\rr$, with $g$ positive definite, $H\not\equiv0$,
and $H(0)=0$. Suppose that
\[
 \supp\widehat H\subseteq\{k\in\zz_D:\widehat g(k)=0\}.
\]
Then $g+\varepsilon H$ is not positive definite for any real
$\varepsilon\ne0$.
\end{proposition}
\end{samepage}
\begin{proof}
Suppose that $g+\varepsilon H$ were positive definite. At a Fourier
zero of $g$, this would require $\varepsilon\widehat H(k)$ to be real
and nonnegative. Outside that zero set, $\widehat H$ vanishes by
assumption. Every Fourier value of $\varepsilon H$ would therefore
be real and nonnegative, and inversion at the origin would give
\[
 0=\varepsilon H(0)
   =\frac1D\sum_{k\in\zz_D}\varepsilon\widehat H(k).
\]
The summands are nonnegative and have sum zero, so they must all
vanish. Fourier inversion then gives $H=0$, a contradiction.
\end{proof}

To apply the proposition to Case~2 of
Theorem~\ref{thm:primepreserving}, recall that the perturbation $H$
has pure Fourier order $qN$, while the layered complement $C$ vanishes
at that order. The Fourier transform $\widehat g$ of its normalized
autocorrelation therefore vanishes on $\supp\widehat H$, and $H(0)=0$. The same argument applies
to the perturbation of the period-$180$ example: its Fourier support
has order $36$, at which the autocorrelation term vanishes.

The appendix construction has the same feature in product coordinates.
Under the identification
$\zz_D\simeq\zz_M\times\zz_r$, the support of $\widehat H$ consists
of characters whose old component has order $N$. Each translated
layer of $B^\uparrow$ vanishes at every such old character because
$\Ph_N\mid B$. Hence $\widehat g$ is zero on this entire support,
and again $H(0,0)=0$. Proposition~\ref{prop:pdperturbation} applies.
Changing an ordinary polynomial representative, as in the coefficient
transfer lemmas, leaves these cyclic functions and their Fourier
transforms unchanged.

Symmetrization does not remove this obstruction while retaining a
nonzero perturbation. Indeed, let
$H_{\mathrm{ev}}(x)=(H(x)+H(-x))/2$. Since $g$ is real and positive
definite, it is even and its Fourier zero set is invariant under
negation. Thus $H_{\mathrm{ev}}(0)=0$ and its Fourier support is still
contained in that zero set. If $H_{\mathrm{ev}}\ne0$, the proposition
applies; if $H_{\mathrm{ev}}=0$, symmetrization returns $g$ and restores
the mixed Fourier zero.

The preceding proposition concerns the form of the perturbation.
For the five-point set used in the explicit example, we can make a
stronger statement about the complements themselves. Restriction to
the subgroup corresponding to the original period allows us to use
the rigidity argument from Proposition~\ref{prop:fixedmask}.

\begin{samepage}
\begin{corollary}\label{cor:pdset}
Let $E=\{0,2,3,4,6\}$ and $q\ge1$. Suppose that $P\in\qq[X]$ has
nonnegative coefficients, satisfies $P(0)=1$ and ordinary $(T1)$, and
obeys
\[
 E(X^q)P(X)\equiv U_{30q}(X)\mod (X^{30q}-1).
\]
If its cyclic coefficient function on $\zz_{30q}$ is positive
definite, then $P$ satisfies ordinary $(T2)$. In particular, this
exclusion applies to $A_*=6E\subset\zz_{180}$.
\end{corollary}
\end{samepage}
\begin{proof}
We first restrict the complement to a copy of $\zz_{30}$, where its
values are determined by the set rigidity argument. Set $D=30q$ and
$f=w_P^D$. Nonnegativity and $P(0)=1$ give $f(0)\ge1$, whereas the
convolution equation at zero gives $f(0)\le1$. Define
$h(x)=f(qx)$ on $\zz_{30}$. Restricting the convolution equation to
$q\zz_D$ gives
\[
 h\ge0,\qquad h(0)=1,\qquad
 \one_E*h=\one_{\zz_{30}}.
\]
The set rigidity argument in Proposition~\ref{prop:fixedmask} yields
$h=\one_{5\zz_{30}}$, so $\widehat h$ vanishes at order six.
Restriction to the subgroup also gives the Fourier identity
\begin{equation}\label{eq:pdsetrestriction}
 \widehat h(\ell)=\frac1q
   \sum_{\substack{k\in\zz_D\\k\equiv\ell\bmod 30}}\widehat f(k),
 \qquad \ell\in\zz_{30}.
\end{equation}
Positive definiteness makes every summand nonnegative. Thus a positive
Fourier value of $f$ contributes a positive Fourier value at the order
of its restriction, with no possibility of cancellation.

We use this observation to exclude a failure of $(T2)$. If $P$ failed
that condition, Lemma~\ref{lem:weakallocation} would give an order
$d\mid D$ whose constituent prime powers all belong to $S_P^*$ but
for which $\Ph_d\nmid P$. Choose a character of order $d$ at which
$\widehat f(k)>0$, and put $e=d/\gcd(d,q)$. The tiling identity forces
$E(X^q)$ to vanish at this character, so $E$ vanishes at order $e$.
Since $E(X)=\Ph_5(X)\Ph_6(X)$, we have $e\in\{5,6\}$.

If $e=5$, then $5^{1+v_5(q)}$ is a constituent prime power of $d$ and
belongs to the signature of $E(X^q)$ by
\eqref{eq:dilatedsignature}. This contradicts allocation, since every
constituent of $d$ belongs to $S_P^*$. Therefore $e=6$. The residue
$\ell=k\bmod30$ then has character order six. Equation
\eqref{eq:pdsetrestriction} gives $\widehat h(\ell)>0$, contradicting
the order-six vanishing of $h$. Hence $P$ satisfies $(T2)$.
\end{proof}

The proposition explains the limitation of the signed perturbations,
and the corollary gives an exclusion for the particular set $E$.
Whether other set tiles admit normalized nonnegative positive-definite
$(T1)$ complements failing $(T2)$ remains open. Case~1 of
Theorem~\ref{thm:primepreserving} is separate from this perturbation
argument: it starts from a set complement already failing $(T2)$.
Such a complement would also give a positive-definite weak complement
by normalized autocorrelation, as discussed in
Remark~\ref{rem:newnonpd}.

\begin{samepage}
\section*{Acknowledgements}
G. K. was supported by grants NKFIH-142993, NKFIH-150576, and NKFIH-154121.

C. M. received support for postdoctoral research during the
\textit{Simons Special Semester} at the R{\'e}nyi Center and was also
supported by an Arts \& Science PDF Award and an NSERC Canadian PDF Award.

M. M. was supported by grants NKFIH-146387 and NKFIH-154121.

G. S. was supported by a Discovery Project grant (DP250104965) of the Australian Research Council and grant NKFIH-150576.
\end{samepage}

\appendix
\section{Normalization by shifted layers}\label{sec:lift}
We give an alternative construction on $\zz_M\times\zz_r$, where
$r$ is a prime coprime to $M$. The product coordinates allow us to
place the filter in a nonzero layer, leaving the value at the origin
unchanged. As before, translations of the set complement provide
enough positive weight to absorb the negative filter coefficients.

\subsection{The general construction}
Assume throughout this subsection that $A\oplus B=\zz_M$, that $A$
satisfies $(T2)$, and that
\begin{equation}\label{eq:commonzero}
 N=\prod_{i=1}^k p_i^{\alpha_i}\mid M,\qquad k\ge2,
 \qquad \Ph_N\mid A,\quad \Ph_N\mid B,\quad
 \Ph_{p_i^{\alpha_i}}\mid B\quad(1\le i\le k).
\end{equation}
In particular, the zero of $B$ at the common mixed order $N$ is one of
the zeros required by $(T2)$. We will remove this zero while retaining
the constituent prime-power zeros. The additional cardinality condition
in the definition of an unsupported divisor is not needed for this
construction.

Let $R=R_{N,M}$, let $S_-=\{x:R(x)=-1\}$, and choose a prime $r>2$
coprime to $M$. Put
\[
 D=rM,\qquad K=r|B|,
 \qquad \zz_D\simeq\zz_M\times\zz_r.
\]
Here and below this is the literal residue-coordinate identification:
$z\mapsto(z\bmod M,z\bmod r)$. Choose shifts $t_0,\ldots,t_{r-1}\in\zz_M$
and define
\begin{equation}\label{eq:layeredsets}
 A^\uparrow=A\times\{0\},\qquad
 B^\uparrow=\bigcup_{j\in\zz_r}(B+t_j)\times\{j\}.
\end{equation}
Each layer is a translated copy of the base tiling, and consequently
$A^\uparrow\oplus B^\uparrow=\zz_D$. We now take the normalized
autocorrelation of this complement and place the filter in layer one:
\begin{equation}\label{eq:ghf}
 g=\frac1K\one_{B^\uparrow}*\one_{-B^\uparrow},\qquad
 H(x,j)=R(x)\one_{\{1\}}(j),\qquad f=g+\varepsilon H.
\end{equation}

Let $F$ be the canonical mask polynomial of $f$. The following lower bound on $g$
will ensure nonnegativity; we arrange it in Lemma~\ref{lem:layers}.

\begin{samepage}
\begin{theorem}\label{thm:general}
Under \eqref{eq:commonzero}, suppose that $\varepsilon\in\qq_{>0}$ and
\begin{equation}\label{eq:coverage}
 g(x,1)\ge\varepsilon\qquad(x\in S_-).
\end{equation}
Then the function in \eqref{eq:ghf} satisfies
\[
 f\ge0,\qquad f(0,0)=1,\qquad \one_{A^\uparrow}*f=\one_{\zz_D}.
\]
Its mass is $K$, and its period-relative prime-power zero set is
\[
 S_{F,D}=S_B^*\cup\{r\}.
\]
It satisfies period-relative $(T1)$ and fails $(T2)$ at order $N$.
The set $A^\uparrow$ satisfies ordinary $(T1)$ and $(T2)$, while $f$ is
not positive definite.
\end{theorem}
\end{samepage}
\begin{proof}
By Lemma~\ref{lem:autocorrelation}, $g$ has the required tiling identity,
mass, and value at the origin. The perturbation preserves the tiling
identity because
$\one_{A^\uparrow}*H=(\one_A*R)\otimes\one_{\{1\}}=0$ by
Lemma~\ref{lem:purefilter}. Its negative sites are exactly the points
$(x,1)$ with $x\in S_-$, where \eqref{eq:coverage} ensures $f\ge0$. Finally, $H(0,0)=0$ and $\mass(H)=0$, so the
normalization and mass of $g$ are also preserved.

Each translated layer of $B^\uparrow$ retains the prime-power zeros of
$B$, while equal layer sizes give the zero at order $r$. Their
cyclotomic values at one have product
$r|B|=K=|B^\uparrow|$, by Lemma~\ref{lem:allocation}. Consequently these
are all the prime-power factors of its integer mask, and autocorrelation
preserves exactly these Fourier zeros. The perturbation cannot change
them: a character in the support of $\widehat H$ has old component of
mixed order $N$, and hence cannot have prime-power order. Thus $f$ has precisely the
asserted prime-power zeros, and their values at one give its mass $K$.

At the mixed order $N$, the situation is different. Let $\chi$ be any
base character of order $N$, and take the trivial
character in the new coordinate. Since $\Ph_N\mid B$, each translated
layer has zero Fourier value at $\chi$. Hence $\widehat g(\chi,0)=0$,
whereas
\[
 \widehat f(\chi,0)=\varepsilon\widehat R(\chi)\ne0.
\]
Since the constituent prime-power zeros in \eqref{eq:commonzero} were
preserved, this nonzero value gives the asserted failure of $(T2)$.

For the set factor, every old prime-power zero of $A$ is inherited by
$A^\uparrow$. These factors already account for its unchanged
cardinality, so there are no additional prime-power factors. Every
product required by $(T2)$ is therefore an old order, at which vanishing
is inherited from $A$. To finish, we check that the perturbation destroys
evenness. The autocorrelation $g$ is even and $R(0)=1$, while
$r>2$, so
\begin{equation}\label{eq:asymmetry}
 f(0,1)-f(0,-1)=\varepsilon.
\end{equation}
Thus $f$ is not even and cannot be PD.
\end{proof}

We can arrange the required lower bound by prescribing the
differences between successive layer shifts.

\begin{samepage}
\begin{lemma}\label{lem:layers}
For every prime $r>2$ coprime to $M$ with $r-1\ge|S_-|$, one can
choose shifts such that $g(x,1)\ge |B|/K=1/r$ for every $x\in S_-$.
Consequently, condition~\eqref{eq:coverage} holds whenever
$0<\varepsilon\le1/r$.
\end{lemma}
\end{samepage}
\begin{proof}
List $S_-$ as $x_1,\ldots,x_s$. We prescribe these differences as
successive increments of the shifts, using the last increment to close
the cycle. Define
\[
 d_j=x_j\ (1\le j\le s),\quad d_j=0\ (s<j<r),\qquad
 d_0=-\sum_{j=1}^{r-1}d_j\bmod M,
\]
and put $t_0=0$, $t_j=\sum_{h=1}^jd_h\bmod M$ for $1\le j<r$.
Then the cyclic increment $t_j-t_{j-1}$ equals $d_j$, including at
$j=0$. If an increment has the prescribed value $x$, the $|B|$ ordered pairs
\[
 ((b+t_j,j),(b+t_{j-1},j-1)),\qquad b\in B,
\]
all have difference $(x,1)$. Counting these pairs in the normalized
autocorrelation gives the stated lower bound. Since only finitely many
primes divide $M$, we may choose a prime $r$ as large as required.
\end{proof}

With these shifts, we may in particular take the rational value
$\varepsilon=1/K$. At every negative filter site this gives
$f(x,1)\ge(|B|-1)/K$, while the value at the origin remains unchanged.

\subsection{An ordinary polynomial representative}
To pass from period-relative to ordinary $(T1)$, we use
Lemma~\ref{lem:oddtransfer} to remove cyclotomic factors whose orders
do not divide $D$, without changing the cyclic function.

\begin{samepage}
\begin{corollary}\label{cor:global}
The construction of Theorem~\ref{thm:general} admits an ordinary
nonnegative rational representative $P$ with $P(0)=1$ satisfying ordinary
$(T1)$ and failing $(T2)$, without changing its cyclic coefficient
function.
\end{corollary}
\end{samepage}
\begin{proof}
Since $f(x)\ge\varepsilon>0$ at $x=\CRT(0,1)$, this nonzero residue
can be used for the coefficient transfer. Choose a positive integer $m$
clearing the denominators of $F$ and an odd prime $\ell\mid N$.
Lemma~\ref{lem:oddtransfer} gives
\[
 P(X)=F(X)+\frac{X^{x+D}-X^x}{\ell m}.
\]
This polynomial has the same cyclic coefficient function and constant
coefficient as $F$, and all its cyclotomic factors have orders dividing
$D$. Hence its period-relative $(T1)$ condition is ordinary $(T1)$,
while its failure of $(T2)$ at order $N$ is unchanged.
\end{proof}

\end{document}